\documentclass[12pt,absolute]{mymathart}
\usepackage{mymathmacros}

\usepackage{graphicx,subfigure}

\swapnumbers
\newtheorem{lemma}{Lemma}
\newtheorem{definition}[lemma]{Definition}
\newtheorem{theorem}[lemma]{Theorem}

\newtheorem{corollary}[lemma]{Corollary}

\newtheorem{conjecture}[lemma]{Conjecture}

\newfont{\bbf}{msbm10 at 11pt}

\def\R{\mathbb R}
\def\N{\mathbb N}
\def\C{\mathbb C}
\def\Z{\mathbb Z}
\def\D{\mathbb D}

\newcommand{\disk}{\mathbb {D}}

\newcommand{\diskbar}{\ovl{\disk}}

\newcommand{\Circle}{\mbox{\bbf S}^1}

\newcommand{\ovl}{\overline}

\newcommand{\s}{{\underline s}}
\renewcommand{\t}{\underline t}
\renewcommand{\u}{{\mbox{\tt u}}}

\newfont{\Euler}{eusm10 at 12pt}

\newcommand{\sm}{\setminus}
\newcommand{\hide}[1]{}

\newcommand{\Ray}{G_\s}

\renewcommand{\theta}{\vartheta}

\newcommand{\Stable}{\mathcal{R}}
\newcommand{\B}{\mathcal{B}}
\newcommand{\M}{\mathcal{M}}
\newcommand{\I}{\mathcal{I}}

\newcommand{\Sequb}{\overline{\mathcal{S}}}
\newcommand{\Ctilde}{\tilde{\C}}

\usepackage{hyperref}

\title[Bifurcation Loci]
{Bifurcation Loci of Exponential Maps \\
 and
  Quadratic Polynomials:  \\ Local Connectivity,
Triviality of  Fibers, \\
 and Density of
Hyperbolicity}
\author{Lasse Rempe}
\thanks{We gratefully acknowledge that 
  this work was supported by the European Marie-Curie research training 
  network CODY, the ESF research networking programme HCAA, and
  EPSRC fellowship EP/E052851/1.} 
\author{Dierk Schleicher}

\address{Department of Mathematical Sciences, University of Liverpool, L69 7ZL,
United Kingdom}
\email{l.rempe@liverpool.ac.uk}

\address{School of Engineering and Science,
Jacobs University Bremen (formerly International
University Bremen), Postfach 750 561, D-28725
Bremen, Germany}
\email{dierk@jacobs-university.de}

\begin{document}



\begin{abstract}
We study the bifurcation loci of
quadratic (and unicritical) polynomials and exponential maps.
We outline a proof that the exponential bifurcation locus is connected;
this is an analog to Douady and Hubbard's celebrated theorem that
(the boundary of) the Mandelbrot set is connected.

For these parameter spaces, a fundamental conjecture is that hyperbolic
dynamics is dense. For quadratic polynomials, this would follow from the
famous stronger conjecture that the bifurcation locus (or equivalently the
Mandelbrot set) is locally connected. It turns out that a formally
slightly weaker
statement is sufficient, namely that
every point in the bifurcation locus
is the landing
point of a parameter ray.

For exponential maps, the bifurcation locus is not locally connected.
We describe a different conjecture (triviality of fibers) which naturally
generalizes the role that local connectivity has for quadratic or unicritical
polynomials.
\end{abstract}

\maketitle

\hide{
\section{Introduction}
\label{sec:intro}
}


\section{Bifurcation Loci and Stable Components}
\label{Sec:BifLoci}

The family of \emph{quadratic polynomials} $p_c\colon z\mapsto z^2+c$,
parametrized by $c\in\C$, contains, up to conformal conjugacy, exactly
those polynomials with only a single, simple, critical value (at $c$).
Hence this family is  the simplest
parameter space in the dynamical study of polynomials, and has
correspondingly
received much attention during the last two decades. Similarly,
\emph{exponential maps} $E_c\colon z\mapsto e^z+c$ are, up to
conformal conjugacy, the only transcendental entire functions with only
one singular value (the asymptotic value at $c$).\footnote{%
Often the
parametrization $z\mapsto \lambda e^z$ with $\lambda\in\C^*$ is used
instead. This has the
asymptotic value at $0$ and is conformally equivalent to $E_c$ iff
$\lambda=e^c$. This has the advantage that
two maps $z\mapsto \lambda e^z$ and
$z\mapsto\lambda'e^z$ are conformally conjugate iff $\lambda'=\lambda$.
We prefer the parametrization as $E_c$ not only for the
analogy to the quadratic family, but also because all maps $E_c$ have
the same asymptotics near infinity, and because parameter space is
simply connected, which leads to a more natural combinatorial
description.}
This
simplest transcendental parameter space has likewise been studied
intensively since the 1980s.

In the following, we will treat these parameter spaces in parallel,
unless explicitly stated otherwise; we will write
$f_c$ for $p_c$ or $E_c$.
Following Milnor,
we write
$f_c^{\circ n}$ for the $n$-th iterate of $f_c$. The
map $f_c$  is called \emph{stable} if,
for $c'$ sufficiently close to $c$, the maps $f_c$ and $f_{c'}$ are
topologically conjugate on their Julia sets, and the conjugacy depends
continuously on the parameter $c'$ (the former condition implies
the latter in our setting).
We denote by $\Stable$ the \emph{locus of stability}; that is, the 
(open) set of all 
$c\in\C$ so that $f_c$ is stable.
The set $\Stable$ is open and dense in $\C$ \cite{mss,alexmisha}.

\begin{figure}
\subfigure[$p_c$]{\includegraphics[height=0.32\textheight,trim=00
120 170 120,clip]{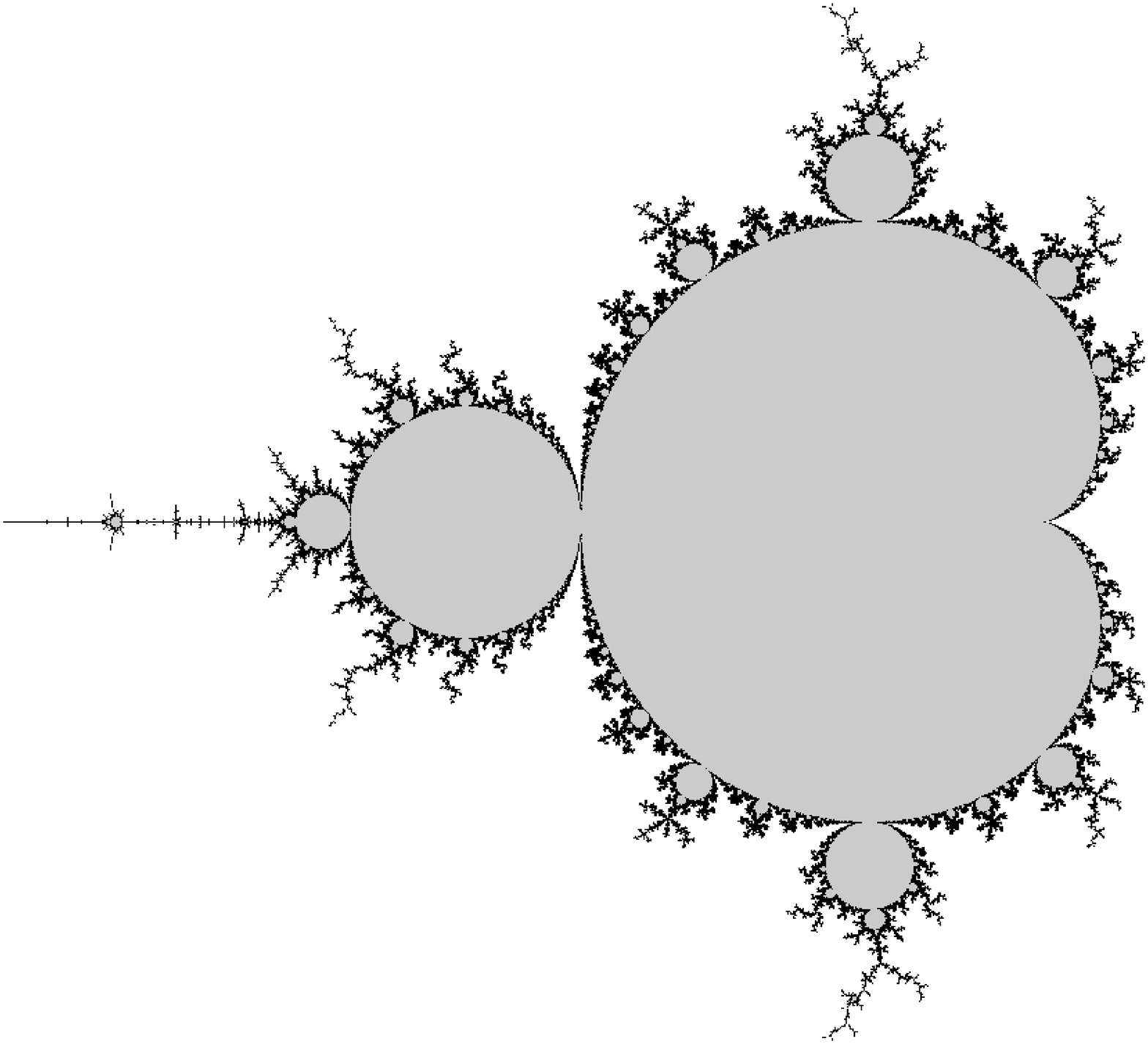}}\hfill
\subfigure[$E_c$]{\includegraphics[height=0.32\textheight]{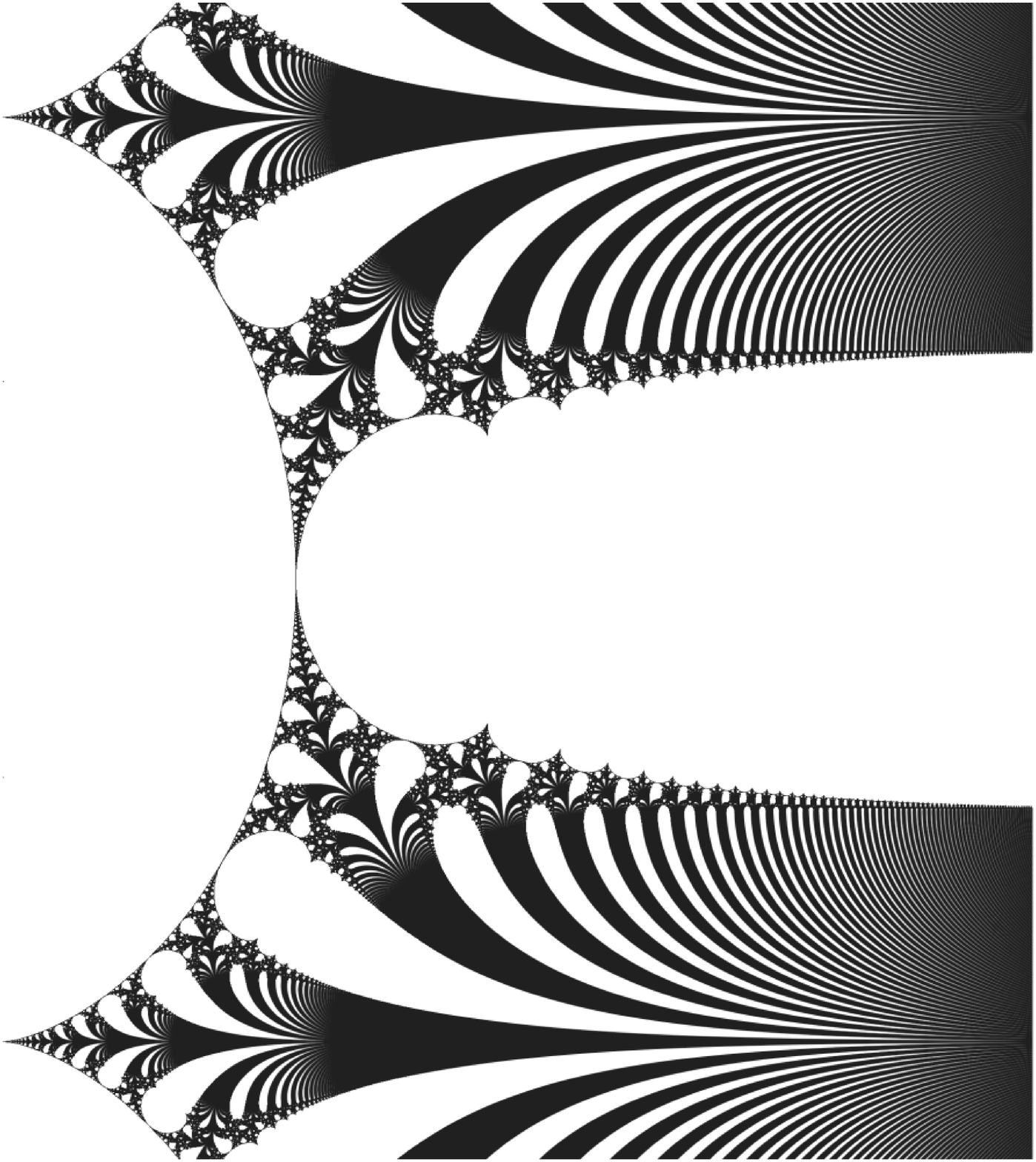}}
\caption{The parameter spaces of quadratic
polynomials $p_c$ (left) and exponential maps
$E_c$ (right); the bifurcation loci are drawn in
black. Unbounded hyperbolic components are white,
bounded hyperbolic components are gray (these do
not exist for exponentials).%
\label{fig:bif}}
\end{figure}

A \emph{hyperbolic component} is a connected component of $\Stable$ in which
every map $f_c$ has an attracting periodic orbit of constant period.
Within any non-hyperbolic component of $\Stable$,
all cycles of $f_c$ would have to be repelling.
One of the fundamental conjectures of one-dimensional
holomorphic dynamics is the following.
\begin{conjecture}[Hyperbolicity is Dense]
\label{Conj:HD}
Every component of $\Stable$ is hyperbolic.
(Equivalently, hyperbolic dynamics is
open and dense in parameter space.)
\end{conjecture}
Hyperbolic components --- both
in the quadratic and in the exponential family ---
are completely understood in terms of their combinatorics
\cite{orsay,expattracting,bifurcations_new}.
The complement $\B:=\C\sm \Stable$ is called the
\emph{bifurcation locus}.  Since $\Stable$ is open and dense in $\C$,
$\B$ is closed and has no interior
points. The bifurcation locus is extremely complicated
(see Figure \ref{fig:bif}).
\begin{theorem}[Bifurcation Loci Connected]
\label{Thm:BifLocusConnected}
$\B$ is a connected subset of $\C$.
\end{theorem}
For quadratic polynomials, $\B$ is the boundary of the famous Mandelbrot
set $\M$, and Theorem \ref{Thm:BifLocusConnected}
is equivalent to the fundamental theorem of Douady
and Hubbard that $\M$ or equivalently
$\partial\M$ is connected \cite[Expos\'e
VIII.I]{orsay}. For
exponential maps, connectivity of the bifurcation locus is new
\cite[Theorem 1.1]{bifurcations_new}; we outline a proof below.

To study the bifurcation locus, it is useful to consider
the \emph{escape locus}
\[
\I:=\{c\in\C\colon f_c^{\circ n}(c)\to\infty \mbox{ as $n\to\infty$}\}
\,\,.
\]
The set $\I$ decomposes naturally into a disjoint
union of \emph{parameter rays}
and their endpoints
(see below),
and the ultimate goal is to describe $\B$ in
terms of these rays. The structure of $\I$ and the
parameter rays is well-understood, and conjecturally,
every point in $\B$ is the landing point of a parameter ray or, in the case
of exponential maps, on a parameter ray itself.
This is in analogy to the dynamical planes of
$f_c$, where the Julia sets are often studied
using the simpler structure of the Fatou set or
of the set of points that converge to $\infty$
under iteration.  We show below 
(Theorem~\ref{Thm:LandingImpliesTrivial}) that,  for the
quadratic family, the conjecture that
every point in $\B$ is the 
landing point of a ray is equivalent to the famous open
question of local connectivity of
the Mandelbrot set; we then reformulate this conjecture in a uniform way for
quadratic polynomials and exponentials.

A fundamental difference between the polynomials $p_c$ and the exponential
maps $E_c$ is their behavior
near $\infty$: every $p_c$ has a superattracting
fixed point at $\infty$ which attracts a neighborhood of $\infty$ in
the Riemann sphere,
while every $E_c$ has an essential singularity at $\infty$ and the set of
points converging to $\infty$ is extremely complicated \cite{expescaping}.
This implies that in the parameter space of
quadratic polynomials we have $\I=\C\sm\M$ and
there is a unique conformal isomorphism
$\Phi\colon \I\to\C\sm\diskbar$ with $\Phi(c)/c\to 1$ as $c\to\infty$
(where $\disk$ is the complex unit disk); the map $\Phi$ was constructed by
Douady and Hubbard~\cite[Expos\'e VIII.I]{orsay} in their proof of
connectivity of the Mandelbrot set. A parameter
ray is defined as the preimage of
a radial line in $\C\sm\diskbar$ under $\Phi$.
On the other hand, the set $\I$ for exponential
maps has a much richer topological structure~\cite{markuslassedierk}:
it is the disjoint  union of uncountably
many curves (\emph{parameter rays}) with or without endpoints; each parameter
ray (possibly with its endpoint) is a path component of $\I$. More 
precisely, every path component of $\I$
is a curve $\Ray\colon(0,\infty)\to \I$ or $\Ray\colon[0,\infty)\to \I$, both
times with $\Ray(t)\to\infty$ as $t\to\infty$.\footnote{%
The parametrization in \cite{markuslassedierk}
is somewhat different, but this is of
no consequence in the following.}
The index $\s$ distinguishes different parameter rays: it is a sequence
of integers and is called the \emph{external address} of $\Ray$.
Different rays have different external addresses, and the set of allowed
external addresses can be described explicitly
\cite{markusdierk,markuslassedierk}.
We call the image $\Ray(0,\infty)$ the \emph{parameter ray} at external
address $\s$, and $\Ray(0)$ its \emph{endpoint} (if it exists). Let $\I_R$
be the union of all parameter rays and $\I_E$ be the set of all endpoints in
$\I$. We say that a parameter ray \emph{lands} if
$\lim_{t\searrow 0}\Ray(t)$ exists (note that many parameter rays
$\Ray$ land at non-escaping parameters, and others might not land at all,
but $\I_E$ consists only of those landing points that are in $\I$).

It is a consequence of the ``$\lambda$-lemma'' from \cite{mss}
that $\B=\partial \I$ both for quadratic polynomials
and for exponential maps; see e.g.\ \cite[Lemma 5.1.5]{thesis}.
In particular, the
bifurcation locus of quadratic polynomials is a compact subset of $\C$,
while for exponential maps it is unbounded.

We will also use the \emph{reduced bifurcation locus}
\[
\B^*:=\B\sm \I_R \,\,;
\]
for quadratic polynomials, clearly $\B^*=\B$, but
for exponential maps,
$\B^*$
is a proper subset of $\B$.

We should remark also on the family of \emph{unicritical} polynomials:
those which have a unique critical point in $\C$.
Such polynomials may be viewed as the topologically and combinatorially
simplest polynomials of a given degree $d$,
and they are affinely conjugate to $p_{d,c}\colon z\mapsto z^d+c$ or to
$z\mapsto \lambda (1+z/d)^d$ with $\lambda=dc^{d-1}$.
These unicritical
polynomials are often viewed as combinatorial interpolation between
quadratic
polynomials and exponential maps; see \cite{dgh,intaddrnew,fibers}.
Everything we say about quadratic polynomials
remains true also
for unicritical polynomials, but for simplicity of exposition
and notation
we usually speak only of quadratic polynomials and of exponentials.

\subsection*{Structure of the article}
In Section \ref{Sec:LocConnTrivFibers}, we review the famous ``MLC''
conjecture for the Mandelbrot set, and then define ``fibers''
(introduced for the case of $\M$ in \cite{fibers})
of quadratic polynomials and
exponential maps in parallel. This concept allows us to formulate
Conjecture \ref{Conj:FibersTrivial} on triviality of fibers, which is
equivalent to MLC in the setting of quadratic polynomials.
We also discuss a number of basic results on fibers. For ease
of exposition, the theorems stated in Section
\ref{Sec:LocConnTrivFibers} will be proved,
separately, in Section \ref{Sec:Proofs}.

Apart from a few somewhat subtle topological considerations, the proofs
are not too difficult, but rely on a number of recent non-trivial
results on the structure of exponential parameter space. We cannot
comprehensively review all of these in the present article,
but have attempted
to present the proofs so that they can be followed
without detailed knowledge of these papers.

While most of the results stated are well-known in the case of the
Mandelbrot set, some observations seem to be new even in this
case. (Compare, in particular, Theorem \ref{Thm:LandingImpliesTrivial},
which
allows a simple, and formally weaker, restatement of the MLC
conjecture.)


\section{Local Connectivity and Trivial Fibers}
\label{Sec:LocConnTrivFibers}

\subsection*{Local connectivity of bifurcation loci}
It was conjectured by Douady and Hubbard that the Mandelbrot set $\M$ is
locally connected; this is perhaps \emph{the} central open problem
in holomorphic dynamics. The following is an equivalent formulation.

\begin{conjecture}[MLC]
\label{Conj:MLC}
The quadratic bifurcation locus $\B=\partial \M$ is locally connected.
\end{conjecture}

One of the reasons that this conjecture is important is that it
implies Conjecture~\ref{Conj:HD} for quadratic polynomials:
see Douady and Hubbard~\cite[Expos\'e XXII.4]{orsay}
(see also \cite[Corollary~4.6]{fibers}, as well as Theorems \ref{Thm:Fibers_HD}
and  \ref{Thm:TrivialFibersMLC} below).
\begin{theorem}[MLC Implies Density of Hyperbolicity]
\label{Thm:MLC_HD}
If the Mandelbrot set is locally connected, then
hyperbolic dynamics is dense in the space of quadratic polynomials.
\end{theorem}

In topological terms, the situation in exponential parameter space is
very different from what we expect in the Mandelbrot set: the analog of
Conjecture \ref{Conj:MLC} is known to be false.

\begin{theorem}[Failing Local Connectivity of Exponential Bifurcation Locus]
\label{Thm:ExpoNonLocConn}
The exponential bifurcation locus $\B$ is not locally connected.
More precisely, $\B$ is not locally connected at any point of $\I_R$.
\end{theorem}

\begin{figure}
\begin{center}
\includegraphics[width=\textwidth]{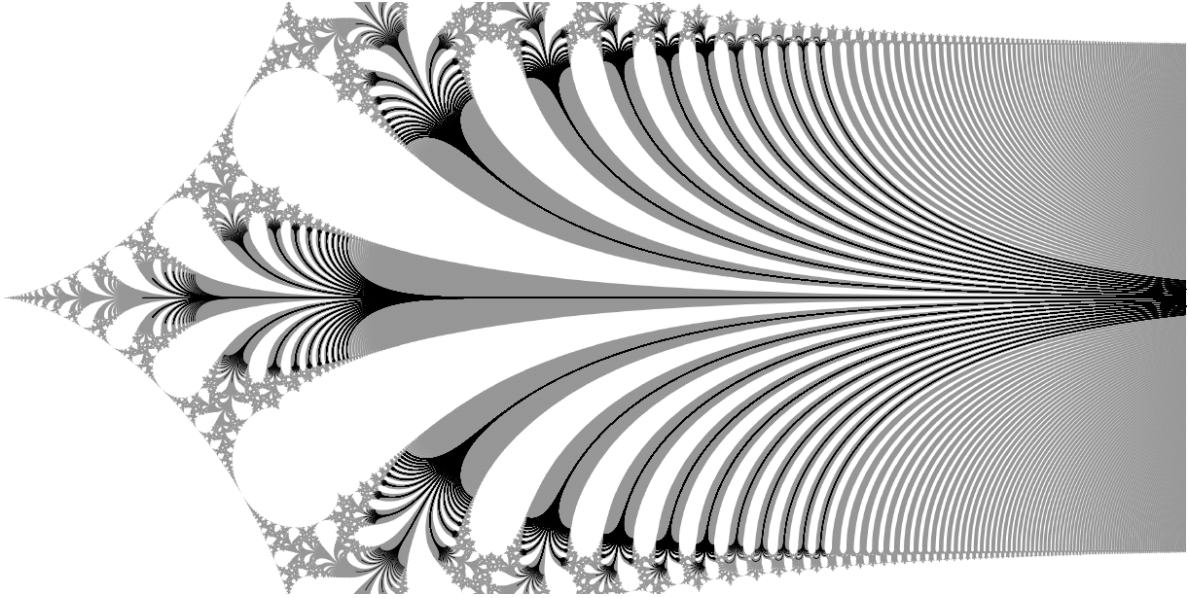}
\end{center}
\caption{Exponential parameter space contains
``Cantor bouquets'', which are closed sets consisting of
uncountably many disjoint simple curves. Some such curves are
indicated here in black.
\label{fig:bouquet}}
\end{figure}

In essence, failure of local connectivity of $\B$ in the exponential
setting was first shown by Devaney:
from his proof  \cite{devaneyunstable}
that the exponential map $\exp$
itself is not structurally stable, it also follows that $\B$ is not
locally connected at the point $c=0$.
Theorem \ref{Thm:ExpoNonLocConn}
is related to the existence of so-called
\emph{Cantor bouquets} within the bifurcation locus; compare Figure
\ref{fig:bouquet}. These consist of
uncountably many curves (on parameter rays) that are locally modelled as
a subset of the product of an interval and a Cantor set.

\subsection*{Fibers}
Local connectivity of the Mandelbrot set would have a number of
important consequences apart from that described by
Theorem \ref{Thm:MLC_HD}: for example,
there are a number of topological models for $\M$
(such as Douady's \emph{pinched disks}
\cite{pincheddisk} or Thurston's \emph{quadratic
minor lamination} \cite{thurstonlaminations})
that are homeomorphic to $\M$ if and only if $\M$
is locally connected. Moreover,
MLC is equivalent to the \emph{combinatorial
rigidity conjecture}: any two maps in $\B$
without indifferent periodic orbits
can be distinguished combinatorially (in terms of
which periodic dynamic rays land together).
Hence it is desirable to find another topological
concept which can play the role of
local connectivity in the space of exponential maps. One convenient
notion of this type is
provided by \emph{triviality of fibers}, introduced for the
Mandelbrot set in \cite{fibers}; it has the advantage of transferring
easily to exponential parameter space.
Another advantage is that even for polynomials,
any possible failure of local connectivity
can be made more precise by giving topological
descriptions of non-trivial fibers.

\begin{definition}[Separation Line and Fiber]
\label{Def:SeparationLine}
A \emph{separation line} is a Jordan arc $\gamma\subset\C$ in parameter
space, tending to $\infty$ in both directions and containing
only hyperbolic and finitely many parabolic
parameters\footnote{Since all hyperbolic
components of exponential parameter space are
unbounded and $\infty$ is accessible,
it suffices to allow just a single parabolic
parameter on every separation line; for $\M$ at
most
two parabolic parameters suffice.}.

A separation line $\gamma$ \emph{separates}
two points $c,c'\in \B^*$ if
$c$ and $c'$ are in two different components of $\C\sm\gamma$.

The \emph{extended fiber} of a point $c\in\C$ is the set of all $c'$ which
cannot be separated from $c$ by any separation line. The
\emph{(reduced) fiber} of $c\in\C\setminus\I_R$
consists of all $c'$ in the extended fiber of $c$ which do not
belong to $\I_R$.

A fiber is called \emph{trivial} if it consists of exactly one point.
\end{definition}
\begin{remark}[Remark 1]
This definition is somewhat different from (but equivalent to)
that originally given in
\cite{fibers} for quadratic polynomials. See the remark on alternative
definitions of separation lines below.
\end{remark}
\begin{remark}[Remark 2]
One fundamental difference between polynomial and
exponential parameter space is that the escape
locus $\I$ is hyperbolic only in the polynomial
case; for exponential parameter space, separation
lines are disjoint from $\I$ and thus from
parameter rays (an alternative definition of
separation lines uses parameter rays; see below).
\end{remark}
\begin{remark}[Remark 3]
   If $c$ is a hyperbolic parameter, then it follows easily from the
    definition that the fiber of $c$ is trivial. Hence all interest lies
    in studying non-hyperbolic fibers, and we will usually restrict our
    attention to this case.
\end{remark}

\begin{figure}
\includegraphics[width=\textwidth]{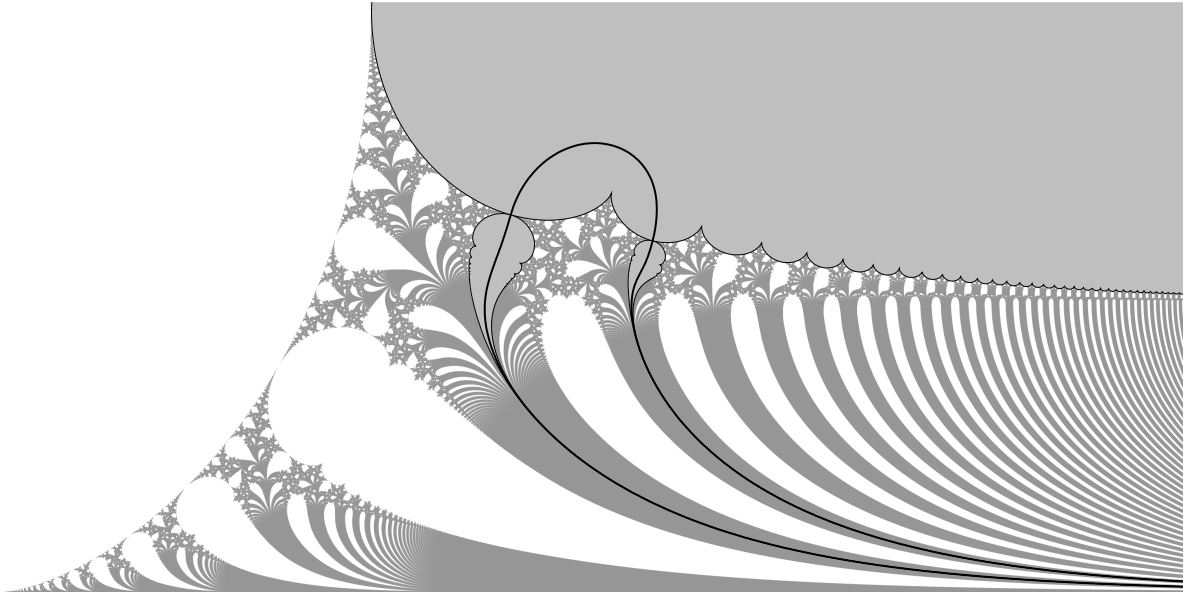}
\caption{A separation line in the space of exponential maps}
\end{figure}

\begin{theorem}[Properties of fibers] \label{thm:fiberproperties}
Every extended fiber $\breve{Y}$
  is a closed and connected subset of parameter space. 
Either $\breve 
Y$ is a single point in a hyperbolic component, or
  $\breve{Y}\cap\B^*\neq\emptyset$.

  In particular, $\breve{Y}$
   contains the (reduced) fiber $Y = \breve{Y}\setminus \I_R$; this fiber
   is either trivial or uncountable.
\end{theorem}
\begin{remark}[Remark 1]
For fibers in the Mandelbrot set, these claims are immediate from
  the definition. For exponential maps, the proofs are more subtle, and
  require some topological considerations as well as rather detailed
  knowledge of parameter space.
\end{remark}
\begin{remark}[Remark 2]
Once we know that
parabolic parameters have trivial fibers
(which is well-known for quadratic polynomials
\cite{fibers,tanleiparabolics,HubbardYoccoz};
for exponential
maps, a proof will be provided in \cite{rationalfibers}),
it also
follows that fibers are pairwise disjoint.
\end{remark}

Armed with the definition of fibers, we can now propose the following
central conjecture.

\begin{conjecture}[Fibers are trivial]
\label{Conj:FibersTrivial}
For the spaces of quadratic polynomials and exponential maps,
all fibers are trivial.
\end{conjecture}

Since any non-hyperbolic stable component would be contained in a
single fiber, it follows immediately that triviality of fibers
would settle density of hyperbolicity:

\begin{theorem}[Triviality of fibers implies density of hyperbolicity]
\label{Thm:Fibers_HD}
Conjecture \ref{Conj:FibersTrivial} implies Conjecture
\ref{Conj:HD}. \qedoutsideproof
\end{theorem}

\subsection*{Triviality of fibers and local connectivity of the Mandelbrot set}

Following Milnor \cite[Remark after Lemma 17.13]{jackdynamicsthird},
we say that a topological
space $X$ is \emph{locally connected at a point $x\in X$} if $x$ has
a neighborhood base consisting of connected sets.\footnote{%
Sometimes this property is instead referred to as ``connected im kleinen''
(cik) at $x$, and the term ``locally connected at $x$'' is reserved
for what Milnor calls ``openly locally connected at $x$''.}
Triviality of the fiber of a parameter $c$ in the quadratic bifurcation locus
implies local connectivity of the Mandelbrot set at $c$. In fact, all known
proofs of local connectivity of $\M$ at given parameters do so by actually
establishing triviality of fibers.

\begin{figure}
\subfigure[]{\includegraphics[width=.45\textwidth,height=.3\textwidth]{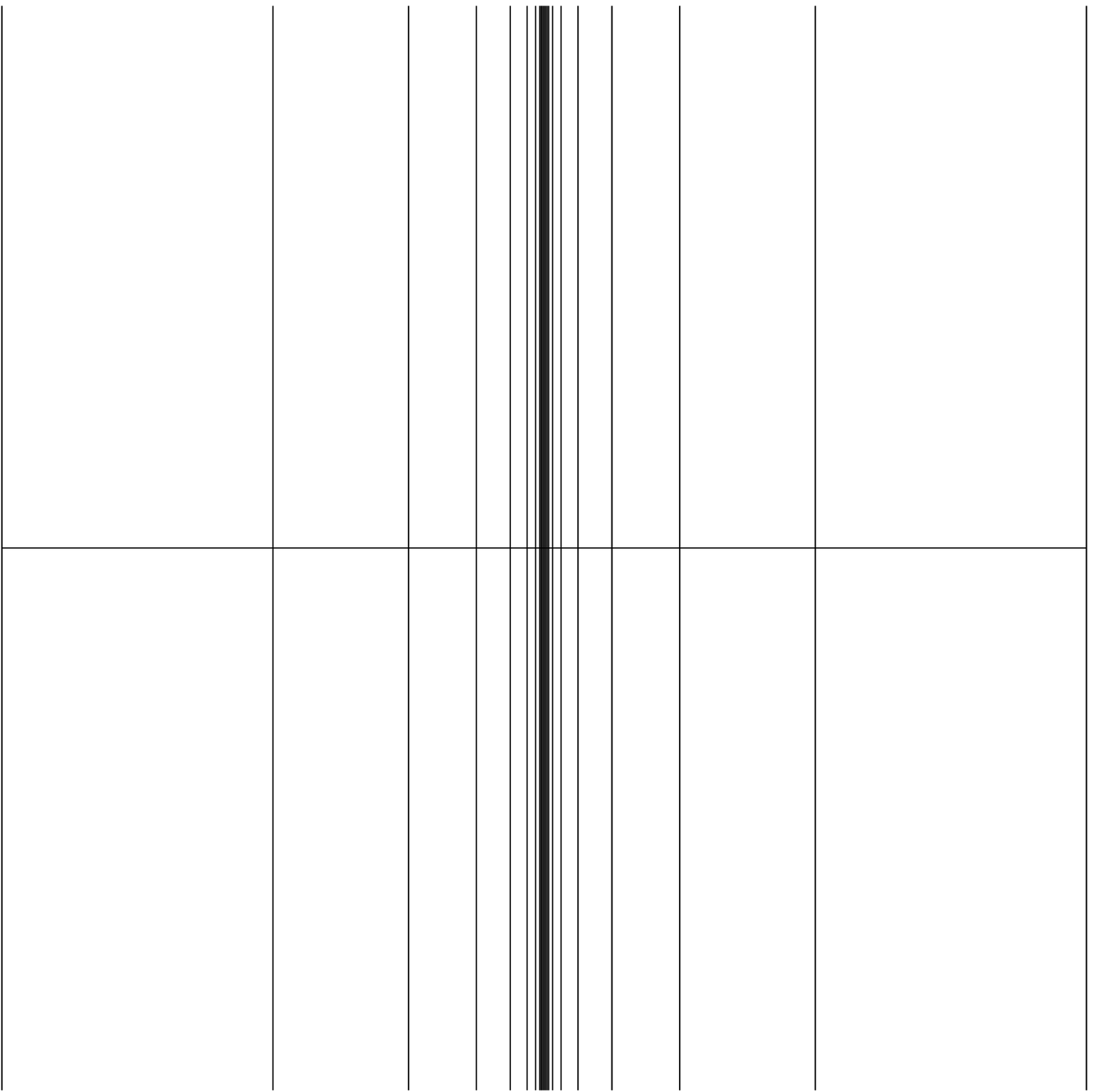}\label{fig:nonaccessible}}\hfill
\subfigure[]{\includegraphics[width=.45\textwidth,height=.3\textwidth]{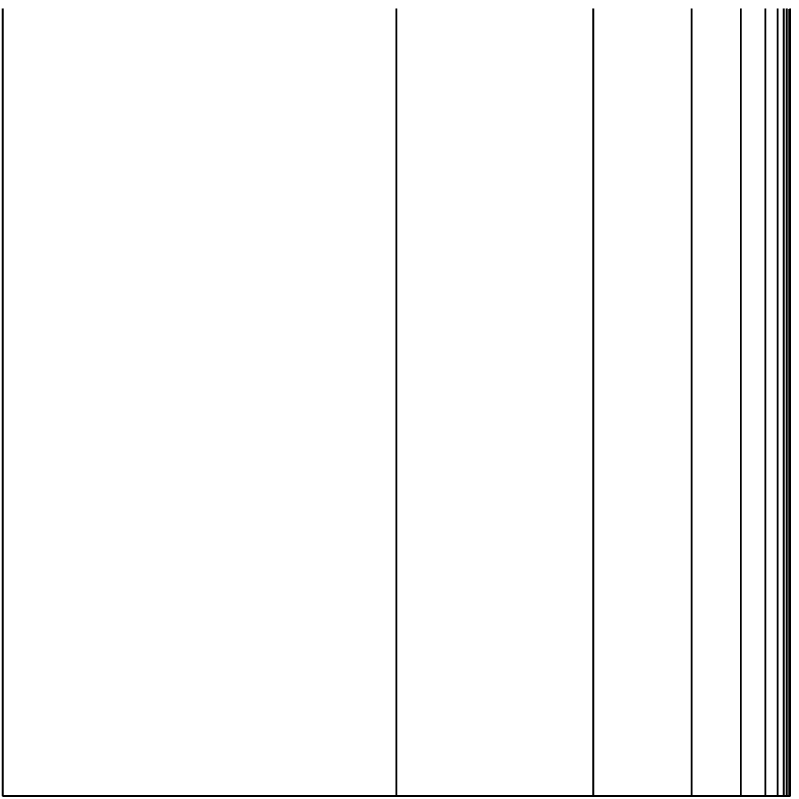}\label{fig:nonlc}}
\caption{Examples of compact connected sets $K\subset\C$ with connected
complement such that
(a) $K$ is locally connected at a point $z_0\in K$, but $z_0$ is not accessible
from $\C\setminus K$
(where $z_0$ is the point at the center of the figure);
(b) every $z\in K$ is accessible from $\C\setminus K$ (and hence the landing
point of an external ray),
but $K$ is not locally connected.}
\end{figure}

The converse question --- in how far local connectivity implies
   triviality of fibers --- is more subtle. For example, a
   full, compact, connected set $K$ may well contain non-accessible points
   at which $K$ is locally connected; see Figure \ref{fig:nonaccessible}.
   On the other hand,
   triviality of the fiber of a parameter $c$ implies that
   there is a parameter ray landing at $c$
   (see Theorem \ref{Thm:LandingImpliesTrivial} below), and hence that $c$
   is accessible from the complement of $\M$. So local connectivity
   at $c$ does not formally imply triviality of the fiber of $c$, but there
   is the following, more subtle, connection.

\begin{theorem}[Trivial Fibers and MLC]
\label{Thm:TrivialFibersMLC}
Let $c$ belong to the quadratic bifurcation locus $\B$. Then the following
are equivalent:
\begin{enumerate}
\item the Mandelbrot set $\M$ is locally connected at every point of the
fiber of $c$;
\item the fiber of $c$ is trivial.
\end{enumerate}

In particular,
local connectivity of the Mandelbrot set is equivalent to
triviality of all fibers in the quadratic bifurcation locus.
\end{theorem}

\subsection*{Triviality of Fibers and Landing of Parameter Rays}

For the Mandelbrot set, by Carath\'eodory's
theorem \cite[Theorem~17.14]{jackdynamicsthird}
local connectivity implies  that every parameter ray
lands, and the map assigning to every external angle the landing point of the
corresponding parameter ray is a continuous surjection $\Circle\to\B$.

Again, replacing local connectivity by triviality of fibers allows us to
obtain a statement regarding the landing of rays
which is true in both quadratic and exponential
parameter space.

\begin{theorem}[Landing of rays implies triviality of fibers]
\label{Thm:LandingImpliesTrivial}
Let $c\in\B^*$. Then the following are equivalent.
\begin{enumerate}
\item The fiber of $c$ is trivial.
\item Every point in the fiber of $c$ is the landing point of a parameter
ray.
\end{enumerate}

In particular, triviality of all fibers is
equivalent to the fact that every point
of $\B^*$ is the landing point of a parameter ray.
\end{theorem}
The last sentence in the theorem provides a convenient way of stating
Conjecture \ref{Conj:FibersTrivial} without the definition of fibers.

The previous result leads to
an interesting observation about the Mandelbrot
set, which is new as far as we know: \emph{MLC is equivalent to the
claim that every $c\in\partial\M$ is accessible
from $\C\setminus\M$}. For general
compact sets, this is far from true (compare Figure
\ref{fig:nonlc}).

\begin{figure}
\begin{center}\resizebox{\textwidth}{!}{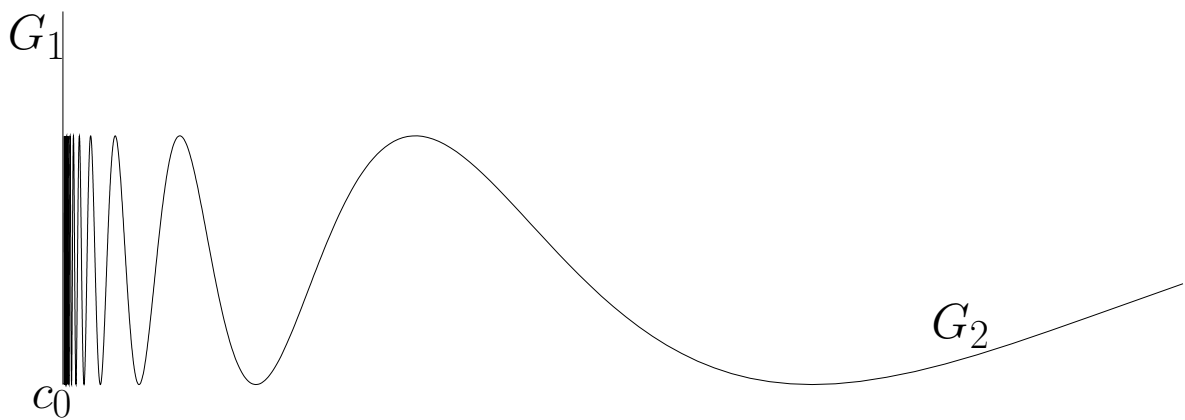}
\end{center}
\caption{It is conceivable that an extended fiber in exponential parameter
space consists of one ray $G_1$ landing at some parameter $c_0$ together
with a second ray $G_2$ which accumulates not only on $c_0$, but also
on a segment of $G_1$. In this case, the fiber of $c_0$ would be trivial,
but the ray $G_2$ clearly does not land at $c_0$.\label{fig:accumulation}}
\end{figure}

In this context, there is
a difference between the situation for
the Mandelbrot set and that of exponential parameter space. In the former
case, triviality of the fiber of $c$ implies that \emph{all} parameter rays
accumulating on $c$ land. In the latter case, we can only conclude
that \emph{one} of these rays lands. The problem is that escaping parameters
are contained in the bifurcation locus;
it is compatible with our current knowledge that
one parameter ray might accumulate on another ray together with its
landing point
(compare Figure \ref{fig:accumulation}), in which case the corresponding
fiber could still be trivial. This might lead us to formulate a stronger
variant of Conjecture \ref{Conj:FibersTrivial}: \emph{all fibers are
trivial, and furthermore all parameter rays land}. It seems plausible that
these conjectures are equivalent (we believe that both are true).

\subsection*{Alternative definitions of separation lines}

We defined separation lines as curves through hyperbolic components.
We did so since this gives a simple definition and does not require
the landing of periodic parameter rays in the exponential family
(a fact that was proved in \cite{habil}, which is however not formally
   published).

There are a number of other
definitions we could have chosen; for example,
\begin{enumerate}
\item
A \emph{separation line} is a curve consisting either of two parameter rays
landing at a common parabolic parameter, or of two parameter
rays landing at distinct parabolic parameters, together with a curve
which connects these two landing points within a single hyperbolic component.
\label{item:separationparabolic}
\item A \emph{separation line} is a curve
as in (\ref{item:separationparabolic}),
except that we also allow two parameter rays landing at a common
parameter for which the singular value is preperiodic.
\item A \emph{separation line} is a Jordan arc, tending to infinity in both
directions, containing only finitely many parameters which are not
escaping or hyperbolic.
\end{enumerate}

These alternative definitions will yield a theory of fibers for which all of
the above results remain true. However, it is a priori conceivable that
separation lines run through parameters with nontrivial fibers, in which
case the fiber of a nearby
parameter $c\in\B^*$ may depend on the definition of
separation lines\footnote{The original definition
in \cite{fibers} for $\M$ is (b), which
is shown to be equivalent to (a); for $\M$, all
parabolic parameters have trivial fibers.}.
On the other hand, the question of
\emph{triviality} of such a fiber, and hence
Conjecture \ref{Conj:FibersTrivial}, is independent of this definition.


\section{Proofs}
\label{Sec:Proofs}

We begin by showing that the exponential bifurcation locus is not
locally connected.

\begin{figure}
\includegraphics[width=\textwidth]{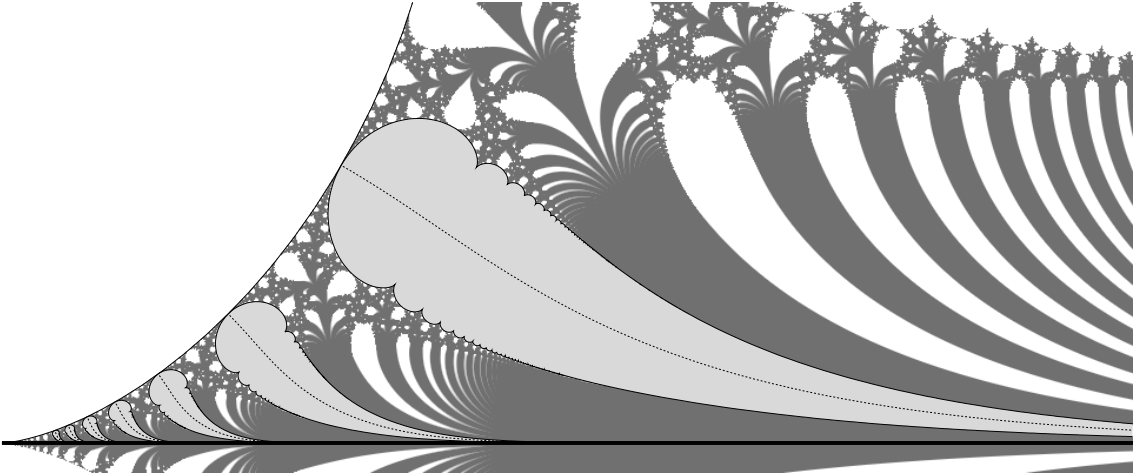}
\caption{Curves in hyperbolic components accumulating on the interval
$[-1,\infty)$ in exponential parameter space.\label{fig:approxreal}}
\end{figure}

\begin{proof}[Proof of Theorem~\ref{Thm:ExpoNonLocConn}]
\emph{(Exponential bifurcation locus not locally connected.)}
As already remarked, failure of local connectivity in the exponential
setting
follows already from Devaney's proof
of structural instability of $\exp$ \cite{devaneyunstable}. Essentially,
he showed
that the interval $[0,\infty)\subset\B$
is accumulated on by curves in hyperbolic
components (compare Figure \ref{fig:approxreal}).

We will now indicate how to prove failure of local connectivity at every
point of $\I_R$ using a similar idea, together with
more detailed knowledge of
exponential parameter space.
More precisely, let
$G_\s\colon\bigl(0,\infty\bigr)\to\C$ be a parameter ray in exponential
parameter space. We claim that, for every $t_0>0$, the curve
$G_\s\colon[t_0,\infty)\to\C$ is the uniform limit of curves
$\gamma_n\colon[t_0,\infty)\to\C$ within hyperbolic components $W_n$ of
period $n$.

Every small neighborhood $U$
of a parameter $G_\s(t)$ with $t>t_0$ then intersects the boundaries of
infinitely many
of these components, and these boundaries are separated from
$G_{\s}(t)$ within $U$
by the curves $\gamma_n$. This proves the theorem.

In the case where $t_0$ is sufficiently large,
the existence of the curves $\gamma_n$
is worked out in detail in \cite[Section 4]{bifurcations_new}. Here, we will
content ourselves with indicating the overall structure of the proof,
for arbitrary $t_0>0$.

Let $c_0 = G_{\s}(t)$ for $t\geq t_0$,
and let $z_n := E^n(c_0) := E_{c_0}^{\circ n}(c_0)$
denote the singular
orbit of $E_{c_0}$. Since $c_0$ is on a parameter ray,
the real parts of the $z_n$ converge to infinity
like orbits under the (real) exponential function. More precisely,
we have the asymptotics
\[  z_n = F^{\circ n}(T) + 2\pi i s_{n+1} + o(1), \]
where  $F(x) = \exp(x)-1$, $\s=s_1s_2s_3\dots$ and $T$ is some positive 
real number.
Given the expansion of $\exp$ in the right half
plane, it should not be surprising that
\[ \left|\left(E^n\right)'( c_0 )\right| \to \infty \]
as $n\to\infty$ (for a proof, see \cite[Lemma
6]{parameterdimensionparadox}). Furthermore, this
growth of the derivative is uniform for $t\geq t_0$.

It follows readily (compare \cite[Lemma 3]{parameterdimensionparadox}) that,
for sufficiently large $n$,
$c_0=G_{\s}(t)$ can be perturbed
to a point $c_n = \gamma_n(t)$ whose singular orbit
follows that of
$E_c$ closely until the $(n-2)$th
iteration, where the two orbits differ by
$i\pi$. In other words,
\[
\left|E^k(c_n)- z_k \right|\ll 1 \quad
\mbox{for all $k=0,1,\dots,n-3$,}
\]
while
\begin{equation}
E^{n-2}(c_n)  =  z_{n-2} +i\pi. \tag{*}
\end{equation}

This implies that the real part of $E_{c_n}^{\circ n-1}(c_n)$
is very negative, and hence $E_{c_n}^{\circ n}(c_n)$ is very close to $c_0$.
So the orbit of the singular value $c_n$ is
``almost'' periodic of period $n$. The contraction
along this orbit is such that a certain disk around $c_n$ is mapped into
itself, and $E_{c_n}$ has an attracting periodic orbit of period $n$
(similarly as in \cite[Lemma 7.1]{bakerexp}; compare
\cite[Lemma 3.4]{expattracting}).
It is easy to see that $c_n=\gamma_n(t)$
is continuous in $t$ and converges uniformly to
$G_\s(t)$ for $n\to\infty$, as desired.
\end{proof}
\begin{remark}[Remark 1]
We remark that the curves $\gamma_n$ constructed in the proof converge
to $G_{\s}\bigl([t_0,\infty)\bigr)$ ``from above'', in the sense that
they tend to infinity in the unique component of
\[ \{\re z > \re G_{\s}(t_0)\} \setminus G_{\s}\bigl([t_0,\infty)\bigr) \]
which contains points with arbitrarily large imaginary parts.
In equation (*), we could have just as well chosen $c_n$ such that
$E^{n-2}(c_n) = z_{n-2} - i\pi$; in this case the curves $\gamma_n$
would converge to $G_{\s}\bigl([t_0,\infty)\bigr)$ from below. We will use this
fact in the proof of Lemma \ref{lem:accumulation} below.
\end{remark}
\begin{remark}[Remark 2]
We do not currently know anything about local connectivity of $\B$
at points of $\B^*$. It seems conceivable that $\B$ is locally
connected exactly at the points of $\B^*$; but in our view the more relevant
question is whether all fibers of points in $\B^*$ are trivial.
\end{remark}

To begin our discussion of fibers, we note the following elementary
consequences of their definition.
\begin{lemma}[Extended Fibers] \label{lem:basicproperties}
Every extended fiber $\breve{Y}$ is a closed subset of $\C$;
the closure $\hat{Y}$
of $\breve{Y}$ in the Riemann sphere is compact and connected.

(In the space of quadratic polynomials, every extended fiber is boun\-ded
and hence $\hat{Y}=\breve{Y}$;
in the space of exponential maps, every non-hyperbolic
extended fiber is unbounded and hence has
$\hat{Y}=\breve{Y}\cup\{\infty\}$.)
\end{lemma}
\begin{proof}
Let $\breve{Y}$ be the extended fiber of a point $c_0$.
The set of points not separated from $c_0$ by a given separation line
is a closed subset of $\C$. So extended fibers are defined as an
intersection of a collection of closed subsets of $\C$, and therefore
closed themselves.

It is easy to verify that a point $c$ which is separated from $c_0$
by a collection of separation lines can also be separated from
$c$ by a single separation line. Furthermore, one can verify that
only countably many separation lines $\gamma$ are necessary to separate
$c_0$ from all points outside $\breve{Y}$.
For example, we can
require that the intersection of $\gamma$ with any
hyperbolic component $W$ is a hyperbolic geodesic of $W$. Since the
set of parabolic parameters is countable and every separation line runs through
hyperbolic components and finitely many parabolic parameters,
there are only countably many such separation lines.

It follows that $\breve{Y}$ can be written as a nested
countable intersection
of closed, connected subsets of the plane; this proves connectivity.
\end{proof}

Before we prove the remaining
theorems, we require some preliminary combinatorial and topological
considerations in exponential parameter space. These become necessary
mainly because we need to take into account the
possibility of parameter rays accumulating on points of
$\I_R$.

We begin by noting that there is a natural combinatorial compactification of
exponential parameter space, as follows.
   The external addresses considered so far are
    elements of $\Z^\N$; we will sometimes call these
    \emph{infinite} external addresses.
    We also introduce
    \emph{intermediate} external addresses: these have the form
$s_1 \dots s_{n-1} \infty$, where $n\geq 1$,
    $s_1,\dots,s_{n-2}\in\Z$ and $s_{n-1}\in \Z+1/2$;
note that there is a unique intermediate external address of length $n=1$,
   namely $\infty$.
The set of all
   infinite and intermediate external addresses will be denoted $\Sequb$.
The lexicographic order induces
a complete total order on $\Sequb\sm\{\infty\}$
and a complete cyclic order on $\Sequb$;
compare \cite[Section 2]{expcombinatorics}.

We can then
define a natural topology on
$\Ctilde = \C\cup\Sequb$, which has the property
that $G_{\s}(t)\to \s$ in this topology
as $t\to+\infty$. The space $\Ctilde$ is homeomorphic to the closed
unit disk $\overline{\D}$, where $\C$ corresponds to the interior of
the disk, and $\Sequb$ to the unit circle.
(This construction is analogous to
compactifying parameter space of quadratic polynomials by adding a circle
of external angles
at infinity.) See also \cite[Appendix A]{expcombinatorics}.

For any parameter ray $G_{\s}$, we now set
$G_{\s}(\infty) := \s\in\Ctilde$, giving a parametrization
$G_{\s}:(0,\infty]\to\Ctilde$.
   Similarly, recall the curves $\gamma_n$ accumulating on the parameter ray
   $G_{\s}$ as constructed in the proof of 
Theorem~\ref{Thm:ExpoNonLocConn}.
They can be extended continuously 
by setting
   $\gamma_n(\infty)=\s_n$, where $\s_n$ is an intermediate external address
   and $\s_n\to\s$.
   Armed with this terminology,
   we can state and prove the following
   key fact, which we will use repeatedly.

\begin{lemma}[Accumulation on parameter rays] \label{lem:accumulation}
Let $A\subset\Ctilde$ be connected, and let $\tilde{A}$ be the closure of
   $A$ in $\Ctilde$. Suppose that $A$ intersects at most finitely many
   hyperbolic components and that
   $A\cap \Sequb$ contains at most finitely many intermediate
   external addresses.

Suppose that $G$ is a parameter ray, and that
   $G(t_0)\in \tilde{A}$ for some $t_0\in(0,\infty]$. Then either
   $\tilde{A}\subset G\bigl([t,\infty]\bigr)$ for some $0<t\leq t_0$, or
   $G\bigl((0,t_0]\bigr)\subset \tilde{A}$.
\end{lemma}
\begin{proof}
If $G((0,t_0])\subset\tilde A$, then nothing is
left to prove, so we may assume that
there is a $t\in(0,t_0)$ with $G(t)\notin \tilde{A}$ and then show that
$\tilde A\subset G([t,\infty])$.
There is
some $\eps>0$ with $\D_{\eps}(G(t))\cap A = \emptyset$.
Recall from the proof of Theorem \ref{Thm:ExpoNonLocConn} and the
subsequent remark that
$G([t,\infty))$ is accumulated on from above resp.\ below by curves
$\gamma_n^+:[t,\infty)\to\C$ and
$\gamma_n^-:[t,\infty)\to\C$, each contained in a
hyperbolic component
of period $n$. Also recall that we can extend these curves continuously by
   setting $\gamma_n^{\pm} := \s_n^{\pm}$ with suitable intermediate external
   addresses $\s_n^{\pm}$. We then have $\s_n^+\searrow \s$ and
   $\s_n^-\nearrow \s$.

By assumption, $A$ intersects at most finitely many
of the curves $\gamma_n^{\pm}$. So for sufficiently large
$n$, $A$ is disjoint from
\[
K_n:=\ovl{\D_\eps(G(t))} \cup\gamma_n^+\bigl([t,\infty]\bigr)
                           \cup\gamma_n^-\bigl([t,\infty]\bigr)
\,\,.
\]
Let $U_n$ be the component of $\Ctilde\sm K_n$ containing $G(T)$ for
   sufficiently large $T$.
   If $\eps$ is sufficiently small, then $G(t_0)\in U_n$,
   hence $A\subset U_n$ and so
\[
A \subset \bigcap_n U_n \subset G_{\s}([t,\infty))
\]
as desired.
(Compare \cite[Corollary 11.4]{topescapingnew} for a similar proof, using
parameter rays accumulating on $G_{\s}$ instead of the curves
$\gamma_n^{\pm}$; recall Figure \ref{fig:bouquet}.)
\end{proof}

Lemma \ref{lem:accumulation} can be applied, in particular, to
  the accumulation sets
  of parameter rays, as stated in the following corollary.
\begin{corollary}[Accumulation sets of parameter rays]
    \label{cor:accumulation}
   Let $G_{\s}$ be a parameter ray, and let $L$ denote the accumulation
    set of $G_{\s}$ in $\Ctilde$. If $G$ is a parameter ray and
    $G(t_0)\in L$ for some $t_0\in (0,\infty]$, then
    $G\bigl((0,t_0]\bigr)\subset L$.

   (In particular, the accumulation set of $G$ is contained in $L$.)
\end{corollary}
\begin{proof}
  Apply the previous lemma to the sets $A = A_t = G_{\s}\bigl((0,t]\bigr)$
   with $t>0$. Since $A$ contains points at arbitrarily small potentials,
   and since parameter rays are pairwise disjoint, the first alternative
   in the conclusion of the lemma cannot hold. Hence (writing
   $\wt{A_t}$ for the closure of $A_t$ in $\Ctilde$ as before) we have
    \[ G\bigl((0,t_0]\bigr) \subset \bigcap_{t>0} \wt{A_t} = L. \qedhere \]
\end{proof}
We also require
the following fundamental theorem from \cite{bifurcations_new}
about exponential parameter space.

\begin{theorem}[The Squeezing Lemma (connected sets version)]
\label{thm:squeezinglemma}
Let $A$ be an unbounded connected subset of exponential parameter space which
contains only finitely many indifferent and no hyperbolic parameters, and let
$\tilde{A}$ be the closure of $A$ in $\Ctilde$.
Then
every $\s\in \wt{A}\cap \Sequb$ is an
infinite external address for which a parameter ray $G_{\s}$
exists, and $G_{\s}(t)\in \wt{A}$ for all sufficiently large $t$.
\end{theorem}
\begin{remark}[Remark 1]
In \cite[Theorem 1.3]{bifurcations_new}, the Squeezing Lemma is
formulated for curves
rather than connected sets. The proof given there also proves the above
version; compare the sketch below.
The name of the result comes from the ``squeezing'' around a parameter ray
by nearby hyperbolic components (or parameter
rays) as in the proof of
Lemma~\ref{lem:accumulation}.

\end{remark}
\begin{remark}[Remark 2]
The idea of the proof of the Squeezing Lemma goes back to the original
proof \cite[Theorem~V.6.5]{habil}, \cite{cras}
that the boundary of every hyperbolic
component
in exponential parameter space is a connected subset of the plane.
A suitable formulation of the Squeezing Lemma can also be used to
prove this fact; see \cite[Theorem 1.2]{bifurcations_new}.
\end{remark}
\begin{proof}[Sketch of proof]
The proof requires a number of technical and combinatorial considerations.
We will describe the underlying strategy and refer the reader to
\cite{bifurcations_new} for details.

Let $\s\in \wt{A}\cap \Sequb$. First, let us show that $\s$
cannot be an intermediate external address. Indeed, otherwise there
is a unique hyperbolic component $W$ which is associated to $\s$
in the following sense    \cite{expattracting}:
if $\gamma:[0,\infty)\to W$ is a curve along which
the multiplier of the attracting orbit tends to zero, then $\s$ is the
sole accumulation point of $\gamma$ in $W$.
Since $A\cap W$ is empty,
we can draw a
separation line through $W$ and two child components (one slightly above
$W$ and one slightly below) which separates $A$ from the address $\s$.
More details can be found in \cite[Section 6]{bifurcations_new}.

So $\s$ is an infinite external address; i.e.\ $\s\in\Z^{\N}$.
If there is no parameter ray associated to $\s$, then
this means \cite{markusdierk} that
at least some of the entries in the sequence $\s$ grow extremely fast.
Every entry in $\s$ which exceeds all previous ones implies the existence of a
separation line in $\Ctilde$ which encloses
parameter rays at external addresses
near $\s$, and the domains thus enclosed shrink to $\s\in\Ctilde$.
Hence
   some of these separation lines separate $A$ from
external addresses near $\s$, a contradiction.
(The detailed proof uses the combinatorial
structure of \emph{internal addresses} and can
be found in \cite[Section 7]{bifurcations_new}.)
It follows that the address $\s$ must indeed have a
parameter ray associated to it.

We can hence apply Lemma \ref{lem:accumulation} with $t_0=\infty$.
    It follows that
    $G_{\s}(t)\in \tilde{A}$ for all sufficiently large $t$, as claimed.
\end{proof}

We show below that every parameter ray has an accumulation point in $\B^*$.
For now, we only note the following.
\begin{lemma}
Every parameter ray has an accumulation point in $\B$. (That is,
a parameter ray cannot land at infinity.)
\end{lemma}
\begin{proof}
For the Mandelbrot set, this is clear. For exponential parameter
space, it follows by applying the Squeezing Lemma to the connected set
$A:=G_{\s}\bigl((0,1]\bigr)$. Indeed, if $A$ is bounded, then
    $G_{\s}$ has at least one finite accumulation point, and there
    is nothing to prove. Otherwise,
    let $\tilde A$ be the closure of $A$ in $\Ctilde$.
    Then there exists an infinite external address $\s'\in \tilde A\cap \Sequb $
    so that $G_{\s'}(t)\in \tilde{A}$ for all sufficiently large $t$. Since
    $A$ does not contain parameters on parameter rays at arbitrarily high
    potentials, it follows that $G_{\s'}(t)$ belongs to
    $\cl{A}\setminus A$, and hence to the accumulation set of $G_{\s}$.
\end{proof}

The following is
essentially a weak version of Douady and Hubbard's ``branch theorem''
\cite[Proposition~XXII.3]{orsay}; see also \cite[Theorem~3.1]{fibers}.

\begin{lemma}   \label{lem:weakbranch}
Let $\breve{Y}$
be a non-hyperbolic
extended fiber. Then $\breve{Y}$ contains the accumulation set of
at least one but at most finitely many parameter rays.
\end{lemma}
\begin{proof}
For quadratic polynomials, the fact that every fiber contains the
accumulation set of a parameter ray is immediate. Indeed, if
$c_0\in\partial\M$, then $c_0$ is contained in some \emph{prime end
impression} (compare e.g.\ \cite[Chapter 17]{jackdynamicsthird} or
\cite[Chapter 2]{pommerenke} for background on the theory of prime ends).
This impression
contains the
accumulation set of an associated parameter ray
(this set is
called the set of \emph{principal points} of the prime end).
Such an accumulation set clearly cannot be separated from $c_0$ by
any separation line.

That there are only finitely many such parameter rays
follows from the usual Branch Theorem, see \cite[Expos\'e XXII]{orsay} or
\cite[Theorem~3.1]{fibers},
which states that the Mandelbrot set
branches only at hyperbolic components and at
postsingularly finite parameters.

For exponential maps, the finiteness statement follows from
the corresponding fact for Multibrot sets. Indeed, suppose we had
   infinitely many parameter rays
   $G_{\s_1}, G_{\s_2},\dots$ which are not separated from one
   another by any separation line. Then it follows by combinatorial
   considerations that all $\s_j$ are bounded sequences, and they differ
   from one another at most by $1$ in each entry.
(The argument is similar as in the proof of the
    Squeezing Lemma: if two addresses
    differ by more than $1$ in one entry, then the ``internal address
    algorithm'' generates a separation line separating the two. Similarly,
    if $\s$ was unbounded, then for every bounded external address
    $\t$ there is
    a separation line which surrounds $\s$ --- i.e., separates $\s$ from
    $\infty$ in $\Ctilde$ --- and separates $\s$ from $\t$. But
    any separation line which surrounds $\s$ must also surround
    all other $\s_j$, which contradicts the fact that there are bounded
    external addresses between any two elements of $\Sequb$.)

So the entries of the addresses $\s_j$ are uniformly bounded,
   and for every Multibrot set $\M_d$ of sufficiently high degree $d$,
   there are
   parameter rays
   $G_{\s_1}, G_{\s_2},\dots$ at corresponding angles. The
   Branch
   Theorem for
   Multibrot sets \cite[Theorem~3.1 and Corollary~8.5]{fibers} shows that
   some of the parameter rays $G_{\s_i}$ are separated
   from each other by a separation line for $\M_d$.
   But then it follows
   combinatorially that there is a similar separation line in exponential
   parameter space, in
   contradiction to our assumption.
   Compare  \cite[Theorem A.3]{expcombinatorics}.

To see that there is at least one parameter ray
contained in every non-hyperbolic
extended fiber $\breve{Y}$, recall from Lemma~\ref{lem:basicproperties} that
$\breve{Y}\cup\{\infty\}$ is a compact and connected
   subset of the Riemann sphere. So every component of
$\breve{Y}$ is unbounded; the Squeezing Lemma
implies that each such component contains
points on some parameter ray $G_{\s}$. But since
separation lines cannot separate the ray $G_\s$,
it follows that $\breve Y$ contains the entire
ray $G_\s$ and hence also its accumulation set.
\end{proof}

\begin{lemma}[Extended fibers are connected]
\label{lem:fiberconnectivity}
Let $\breve{Y}$ be an  extended fiber. Then $\breve{Y}$ and
   $\breve{Y}_{\B} := \breve{Y}\cap\B$ are connected subsets of $\C$.
\end{lemma}
\begin{proof}
We can assume that
$\breve{Y}$ is not a hyperbolic fiber, as otherwise the claim is trivial.

For quadratic polynomials $\breve{Y}$ is a compact and connected subset of $\C$
(Lemma~\ref{lem:basicproperties}). Furthermore, every component of
$\breve{Y}\setminus\B$ (if any) would be a (non-hyperbolic)
stable component of the Mandelbrot set. Since all such components are
simply connected, removing them from $\breve{Y}$ does not disconnect
$\breve{Y}\cap\B$.

It remains to deal with the exponential case.
Let $G_1,\dots,G_n$ be the parameter rays intersecting
   $\breve{Y}$,
with external addresses $\s_1,\dots,\s_n$. Then
all $G_i$ are contained in $\breve Y$.
Let $\tilde{Y}$ be the closure of $\breve{Y}$ in
$\Ctilde$; then, by the Squeezing Lemma,
$\tilde{Y}=\breve{Y}\cup \{\s_1,\dots,\s_n\}$. Note that $\tilde{Y}$
is a compact connected subset of $\Ctilde$ (as in Lemma
\ref{lem:basicproperties}, it can be written as a countable nested
intersection of compact connected subsets of $\Ctilde$).

As in the quadratic setting, any stable component
$U$ of exponential maps is simply connected
   (since boundaries of
    hyperbolic components and parameter rays, both of which are
    unbounded, are dense in $\B$).
Hence $\tilde{Y}_\B := \breve{Y}_{\B}\cup\{\s_1,\dots,\s_n\}$
   is also compact and connected.

Now suppose by contradiction
   that $\breve{Y}_{\B} = A_0 \cup A_1$, where $A_0$ and $A_1$ are
   nonempty, closed and disjoint. Let $\tilde{A}_j$ be the closure of
   $A_j$ in $\Ctilde$. Since $\tilde{Y}_{\B} = \tilde{A}_0\cup\tilde{A}_1$ is
   connected, we must have $\tilde{A}_0\cap \tilde{A}_1\neq \emptyset$.
   I.e., some $\s_j$ belongs to both $\tilde{A}_0$ and $\tilde{A}_1$;
   let $C$ be the component of $\tilde{A}_0$ containing $\s_j$.

By the boundary bumping theorem \cite[Theorem 5.6]{continuumtheory},
   every component of
   $\tilde{A}_0$ contains one of the addresses $\s_i$, so
   $\tilde{A}_0$ has only finitely many connected components.
   This implies that $C\supsetneq\{\s_j\}$. By Lemma \ref{lem:accumulation},
   it follows that $C\cap G_j\neq\emptyset$. Since $G_j$ is connected,
   in fact $G_j\subset A_0$. Likewise, $G_j\subset A_1$,
   which contradicts the assumption that $A_0\cap A_1=\emptyset$.

   So $\breve{Y}_{\B}$ is connected, as is $\breve{Y}$ itself.
\end{proof}

\begin{proof}[Proof of Theorem~\ref{Thm:BifLocusConnected}]
\emph{(Bifurcation locus is connected.)}
Suppose that the bifurcation locus $\B$ is not connected. Then
there is some stable component $U$ such that two components
$C_1$ and $C_2$ of $\partial U$ belong to different connected components
of $\C\setminus U$.
Since boundaries of hyperbolic components are connected subsets of $\C$
\cite{cras,bifurcations_new}
   (recall Remark~2 after Theorem \ref{thm:squeezinglemma}),
$U$ must be a non-hyperbolic stable component, and hence
$\ovl U$ is contained in a single extended fiber $\breve{Y}$
(no separation line can separate any two points in $\ovl U$).
However, then $C_1$ and $C_2$ would belong to different components of
$\breve{Y}\setminus U$, which
contradicts Lemma \ref{lem:fiberconnectivity}.
\end{proof}
\begin{remark}
Theorem~\ref{Thm:BifLocusConnected} can also be proved directly from the
``Squeezing Lemma''; see \cite[Proof of Theorem 1.1]{bifurcations_new}.
\end{remark}

We have now
proved a number of results regarding \emph{extended} fibers. Fibers
themselves can be
more difficult to deal with in the exponential setting, because
they may (at least \emph{a priori}) not be closed. For example, from what
  we have shown so far, it is
  conceivable that an extended fiber is completely
  contained in $\I_R$, and hence does not contain a (reduced) fiber.
  We shall now show
  that this is not the case, and that any extended fiber in fact intersects
  $\B^*$ in either one or uncountably many points.

\begin{theorem}[Accumulation sets of parameter rays]
\label{thm:accumulationsets}
Every parameter ray has at least one accumulation point in $\B^*$;
if there is more than one such point, there are in fact uncountably many.

Furthermore, let
$Y$ be a fiber. Then either
\begin{itemize}
\item $Y$ is trivial (i.e., consists of exactly one point);
if $Y$ is non-hyperbolic, then there is at least one parameter ray
landing at the unique point of $Y$; or
\item $Y$ is not trivial, in which case $Y\cap\B^*$ is uncountable.
\end{itemize}
\end{theorem}
\begin{proof}
For the quadratic family,
the accumulation set of every parameter ray
is contained in the boundary of the Mandelbrot set, which equals $\B^*$.
Any fiber $Y$ is a closed, connected
   subset of $\M$, and hence either consists of a single point or has the
cardinality of the continuum. If the fiber has interior, then its boundary
is contained in $\B^*\cap Y$
and also has the cardinality of the continuum.
If the fiber is trivial, then clearly \emph{every} parameter ray
accumulating on $Y$ must land at the single point of $Y$.

So let us now consider the case of the exponential bifurcation locus.
   Let $G_{\s}$ be a parameter ray, and let
   $\breve{Y}$ be the
   extended fiber containing $G_{\s}$.
   By Lemma \ref{lem:weakbranch}, $\breve{Y}$ contains
   at most
   finitely many parameter rays $G_1,\dots,G_n$, at addresses
   $\{\s_1,\dots,\s_n\}$. (Our original ray $G_{\s}$ will be one of these.)
   Denote by
   $L_i$ the set of all accumulation points of $G_i$ (as $t\to 0$) in
   $\Ctilde$.

   Note that, if $L_i\cap \I_R\neq\emptyset$, say
     $G_j(t) \in L_i$ for some $j\in\{1,\dots,n\}$ and $t\in (0,\infty]$, then
     $G_j\bigl((0,t]\bigr)\subset L_i$ and
     $L_j\subset L_i$ by Corollary \ref{cor:accumulation}.

\textbf{Claim.} \emph{%
   Let $L\subset (\breve{Y}\cap\B)\cup\{\s_1,\dots\s_n\}\subset\Ctilde$
    be nonempty, compact and connected,
    and suppose that there are no $i$ and $t>0$ with
    $L\subset G_{i}\bigl([t,\infty]\bigr)$.}

  \emph{%
   Then $L\cap\B^*$ is either uncountable or a singleton. In the latter case,
    \begin{enumerate}
     \item if $L\not\subset \B^*$, then there is $j$ such that
      $G_j\cap L\neq\emptyset$ and $G_j$ lands at the unique point of
      $L\cap \B^*$; \label{item:raylanding}
     \item every connected
      component of $\C\setminus L$ contains infinitely many parameter rays,
      and hence is not contained in $\breve{Y}$.
       \label{item:separation}
    \end{enumerate}}

   Note that, in (\ref{item:raylanding}), we do not claim that \emph{all}
    rays which intersect $L$ land at the unique point of $L\cap\B^*$.
    To illustrate the claim in this case, it may be useful to
    imagine $L$ to be the set from Figure \ref{fig:accumulation}. Other
    model cases to imagine are those where $L$ is the union of
    of finitely many parameter rays landing at a common point, or where
    $L$ is an indecomposable continuum with one or more parameter rays
    dense in $L$ (as e.g.\ in the standard Knaster
    (or ``buckethandle'') continuum
    \cite[\S 43, V, Example 1]{kuratowski2}).

   Using the claim, we can prove
    both statements of the theorem.
    For the first part, we let $L$ be the accumulation set of 
$G_{\s}$;
by Corollary~\ref{cor:accumulation}, we do not have 
$L\subset G_s\bigl([t,\infty)\bigr)$
for any $i$ and $t>0$, so the 
claim applies.
    For the second part, set $L := (\breve{Y}\cap\B)\cup \{\s_1,\dots,\s_n\}$.
    (The final part of the claim implies, in particular, that if
     $L\cap\B^* = Y\cap\B^*$ is
     a singleton, then $\breve{Y}$ has no interior components, and hence
     $Y$ is trivial, as claimed.)

\begin{subproof}[Proof of the Claim]
   Note that it follows from Lemma \ref{lem:accumulation} and the
    assumption that, for every
    $j\in\{1,\dots,n\}$, there is some $T\in [0,\infty]$ such that
    $L\cap G_j = G_j\bigl((0,T]\bigr)$ (where we understand the interval
     $(0,T]$ to be empty in the case of $T=0$).

   We will proceed by removing isolated end pieces of parameter rays
    from $L$. More precisely, suppose that
    there is some $j\in\{1,\dots,n\}$ and some $t > 0$ such that
    $G_{j}\bigl((t,\infty]\bigr)\cap L$ is a relatively open subset of
    $L$. Choose $t_0\geq 0$ minimal such that all
    $t>t_0$ have this property.

   Then $G_j\bigl((t_0,\infty]\bigr)$ is a relatively
    open subset of $L$, which we will call an ``isolated end piece''.
    Note that the relative boundary of this piece in $L$ is either
    the singleton $G_j(t_0)$, if $t_0>0$, or the accumulation set $L_j$
    of the ray $G_j$ otherwise;
    so in either case this boundary is connected.
    By the boundary bumping
    theorem, $L' := L\setminus G_j\bigl((t_0,\infty]\bigr)$
    is a compact connected
    subset of $\Ctilde$. Furthermore, $L_j$ does not contain any point
    $G_j(t)$ with $t>t_0$ by choice of $t_0$. So
    $L_j\subset L'$, and hence $L'$ is nonempty and
    satisfies the assumptions of the claim. Note that
    $L'\cap \B^* = L\cap \B^*$.

   We can apply this observation repeatedly to remove such isolated end pieces
    of parameter rays from $L$. Note that, if $t_0=0$, then it could be
    that $L'$ contains an isolated end piece of a parameter ray which
    was not isolated in $L$. (Recall Figure
     \ref{fig:accumulation}.) However, since this happens at most $n$
    times, in finitely many steps we
    obtain a set $L_0\subset L$, satisfying the
    assumptions of the claim and with $L_0\cap \B^*=L\cap\B^*$, such that
    furthermore
   \begin{enumerate}
    \item[(*)]
      if $t>0$ is such that
       $G_j\bigl((t,\infty]\bigr)\cap L_0\neq \emptyset$, then
       $G_j\bigl((t,\infty]\bigr)\cap L_0$ is not relatively open in $L_0$.
   \end{enumerate}

   We observe that (*) implies the following stronger property:
    \begin{enumerate}
     \item[(**)]
      if $G_j(t_1)\in L_0$, then $G_j\bigl([t,t_1]\bigr)$ is a nowhere dense
        subset of $L_0$ for all $t\in (0,t_1)$.
    \end{enumerate}
   Equivalently, there are no $t\in (0,t_1)$ and
    $\eps < \min(|t|, |t_1 - t|)$ such that
    $I := G_j((t-\eps,t+\eps))$ is relatively open in $L_0$. To
    prove (**), suppose by contradiction that such $t$ and $\eps$
    exist. By the boundary bumping theorem, every connected
    component of $L_0\setminus I$ must contain one of the two endpoints of $I$.
    Hence there are at most two such components, and these are therefore
    both open and closed in $L_0\setminus I$.
    Furthermore,
    by Lemma \ref{lem:accumulation}, any component of
    $L_0\setminus I$ which intersects $G_j\bigl([t+\eps,\infty]\bigr)$
    is contained in $G_j\bigl([t+\eps,\infty]\bigr)$. Together, these
    facts
    imply that $G_j\bigl((t+\eps,\infty]\bigr)\cap L_0$ is a relatively
    open subset of $L_0\setminus I$, and hence of $L_0$, which contradicts (*).

If $L_0$ is a singleton, then $L_0\subset \B^*$.
   Otherwise, $L_0$ is a nondegenerate continuum, and
   in particular a complete metric space.
   Property (**) implies that $L_0$ can be written as the union of $L_0\cap \B^*$
   with countably many nowhere dense subsets; if $L_0\cap\B^*$ was countable,
   this would violate the Baire category theorem.

In the singleton case, write $L_0=L\cap \B^*=\{c_0\}$.
     Let $I$ be the set of indices $i$ with
     $G_{i}\cap L \neq \emptyset$. Let us assume that $I\neq \emptyset$,
     as otherwise there is nothing to prove. 
     By reordering, we may also assume that
     $I=\{1,\dots,k\}$, where $0\leq k\leq n$,
     and that furthermore $G_1$ is the last
     parameter ray which was completely removed in the construction
     of $L_0$, $G_2$ is the one completely removed before that, etc.
     By construction, the accumulation set $L_1$ of $G_1$ is contained
     in $L_0$, and hence $G_1$ lands at $c_0$.
     In fact, we inductively get
      \begin{equation}
         L_{i} \subset \{c_0\} \cup \bigcup_{j=1}^{i-1} G_{j}.
          \label{eqn:constructionofL}
      \end{equation}
    It follows that every component of
     $\Ctilde\setminus L$ contains an interval of $\Sequb\setminus I$, and
     hence infinitely many parameter rays, as claimed.

    (One way of seeing this is to recall that $\Ctilde$ is homeomorphic to
  the unit disk in $\R^2$. It follows from (\ref{eqn:constructionofL})
     and Janiszewski's theorem (see \cite[Page 2]{pommerenke} or
     \cite[Theorem V.9$\cdot$1$\cdot$2]{newmanplanetopology})
     that $L$, considered as a subset of $\R^2$ in this
     manner, does not separate the plane. So every component of
     $\Ctilde\setminus L$ must intersect the boundary of $\Ctilde$ in $\R^2$,
     i.e. $\Sequb$, as required.)
\end{subproof}
\end{proof}

\begin{proof}[Proof of Theorem \ref{thm:fiberproperties}]
\emph{(Properties of fibers.)}
We just proved
  the fact that fibers are
  either trivial or uncountable.
  Also, we proved that every parameter ray has some accumulation
  point in $\B^*$, and hence that every extended fiber intersects
  $\B^*$.
The fact that extended fibers are
  connected was shown above in Lemma \ref{lem:fiberconnectivity}.
\end{proof}

\begin{proof}[Proof of Theorem~\ref{Thm:LandingImpliesTrivial}]
\emph{(Trivial fibers and landing of rays.)}
Let $Y$ be a fiber, and suppose that every point of $Y$ is the landing
point of a parameter ray. By Lemma \ref{lem:weakbranch},
this means that $Y$ is finite. Hence, by
Theorem \ref{thm:accumulationsets}, $Y$ is trivial.

The converse follows directly from Theorem \ref{thm:accumulationsets}.
\end{proof}

Finally, let us prove the two remaining theorems, which deal
exclusively with the Mandelbrot set $\M$.

\begin{proof}[Proof of Theorem \ref{Thm:TrivialFibersMLC}]
\emph{(Trivial fibers and local connectivity.)}
It is easy to see that triviality of a fiber $Y$ in the Mandelbrot set
implies local connectivity of $\M$ at $Y$. Indeed, as noted above,
$Y$ can be written as the nested intersection of countably many
connected closed subsets of $\B$, each of which is a neighborhood
of $Y$. (Compare \cite[{Proposition~4.5}]{fibers}.)

For the converse direction, suppose that $\M$ is locally connected
at every point of $Y$. Let $z_0 \in Y$ be an accumulation point of
some parameter ray $G$. Then there is a sequence $C_k$
of cross-cuts of
the domain $W := \Ch \setminus\M$ (i.e. $C_k$ is a Jordan arc intersecting
$\M$ only in its two endpoints) with the following properties.
\begin{itemize}
\item $C_k$ separates $\infty$ from all points on $G$ with sufficiently
small potential.
\item The arcs $C_k$ converge to $\{z_0\}$ in the Hausdorff distance.
\end{itemize}

Let $W_k$ be the component of $W\setminus C_k$ not containing $\infty$.
Then $I_G := \bigcap_k \cl{W_k}$ is the \emph{prime end impression}
of the parameter ray $G$. We will show that $I=\{z_0\}$.

Indeed, let $\eps>0$. Since $\M$ is locally connected at $z_0$,
we can find a connected neighborhood $K$ of $z_0$ in $\M$ of diameter less
than $\eps$, for any $\eps>0$. Since closures of connected sets are
connected, we may assume that $K$ is closed.
For sufficiently large $k$, the arc
$C_k$ is a crosscut of $K$, and is contained in the disk of radius
$\eps$ around $z_0$. It follows that $\diam W_k \leq \eps$, and hence
$\diam I \leq \eps$. Since $\eps>0$ was arbitrary, this implies that
$I=\{z_0\}$, as claimed.%
\footnote{In the terminology of prime ends, we have just shown that
the complement of a simply connected domain cannot be
locally connected at any principal point of a nontrivial
prime end impression. Even more is true:
a prime end impression can contain at most two points in which
the complement of the domain
is locally connected; compare \cite{localconnectivity}.}

The boundary of any fiber $Y$ is contained in the union of the prime 
end impressions
corresponding to the parameter rays accumulating at $Y$. By Lemma 
\ref{lem:weakbranch},
   there are only finitely many
   such parameter rays. As we have just shown,
   for each of these rays the prime end impressions consist of a single
    point.  The boundary of fiber $Y$ is thus finite. 
Since $Y$ is 
connected, it follows that $Y$ is trivial as claimed.
\end{proof}

\begin{proof}[Proof of Theorem \ref{Thm:MLC_HD}]
\emph{(MLC implies density of hyperbolicity.)}
By Theorem \ref{Thm:TrivialFibersMLC},
local connectivity of the Mandelbrot set is
equivalent to triviality of fibers; by Theorem \ref{Thm:Fibers_HD},
triviality
of fibers implies density of hyperbolicity.
\end{proof}

\bibliographystyle{hamsplain}
\small{\bibliography{/Latex/Biblio/biblio}}

\begin{thebibliography}{BDG2}

\bibitem[BBS]{parameterdimensionparadox}
Mihai Bailesteanu, H.~Vlad Balan, and Dierk Schleicher, \emph{Hausdorff
  dimension of exponential parameter rays and their endpoints}, 
  Nonlinearity \textbf{21} (2008), 13--20, 
  \mbox{\href{http://www.arXiv.org/abs/0704.3087}{arXiv:0704.3087}}.

\bibitem[BR]{bakerexp}
I.~Noel Baker and Philip~J. Rippon, \emph{Iteration of exponential functions},
   Ann. Acad. Sci. Fenn. Ser. A I Math. \textbf{9} (1984), 49--77.

\bibitem[BDG1]{dghnew1}
Clara Bodel{\'o}n, Robert~L. Devaney, Michael Hayes, Gareth Roberts, Lisa~R.
   Goldberg, and John~H. Hubbard,
   \emph{\href{http://math.bu.edu/people/bob/papers/hairs.ps}{Hairs for the
   complex exponential family}}, Internat. J. Bifur. Chaos Appl. Sci. Engrg.
   \textbf{9} (1999), no.~8, 1517--1534.

\bibitem[BDG2]{dghnew2}
\bysame, \emph{\href{http://math.bu.edu/people/bob/papers/hairs-2.ps}{Dynamical
   convergence of polynomials to the exponential}}, J. Differ. Equations Appl.
   \textbf{6} (2000), no.~3, 275--307.

\bibitem[De]{devaneyunstable}
Robert~L. Devaney, \emph{Structural instability of {${\rm exp}(z)$}}, Proc.
   Amer. Math. Soc. \textbf{94} (1985), no.~3, 545--548.

\bibitem[DGH]{dgh}
Robert~L. Devaney, Lisa~R. Goldberg, and John~H. Hubbard, \emph{A dynamical
   approximation to the exponential map by polynomials}, Preprint, MSRI
   Berkeley, 1986, published as \cite{dghnew1,dghnew2}.

\bibitem[Do]{pincheddisk}
Adrien Douady, \emph{Descriptions of compact sets in {${\bf C}$}}, Topological
   methods in modern mathematics (Stony Brook, NY, 1991), Publish or Perish,
   Houston, TX, 1993, pp.~429--465.

\bibitem[DH]{orsay}
Adrien Douady and John Hubbard, \emph{Etude dynamique des polyn{\^o}mes
   complexes}, Pr{\'e}publications math{\'e}mathiques d'Orsay (1984 / 1985),
   no.~2/4.

\bibitem[EL]{alexmisha}
Alexandre~{\`E}. Eremenko and Mikhail~Yu. Lyubich, \emph{Dynamical properties
   of some classes of entire functions}, Ann. Inst. Fourier (Grenoble)
   \textbf{42} (1992), no.~4, 989--1020.

\bibitem[FRS]{markuslassedierk}
Markus F\"orster, Lasse Rempe, and Dierk Schleicher, \emph{Classification of
   escaping exponential maps}, Proc. Amer. Math. Soc. 136 (2008), 651 -- 663, 
\mbox{\href{http://www.arXiv.org/abs/math.DS/0311427}{arXiv:math.DS/0311427}}.

\bibitem[FS]{markusdierk}
Markus F\"orster and Dierk Schleicher, \emph{Parameter rays for the exponential
   family}, Preprint, 2005,
   \mbox{\href{http://www.arXiv.org/abs/math.DS/0505097}{arXiv:math.DS/0505097}%
}, to appear in Ergodic Theory Dynam. Systems.

\bibitem[H]{HubbardYoccoz}
J.~H. Hubbard, \emph{Local connectivity of {J}ulia sets and bifurcation loci:
   three theorems of {J}.-{C}. {Y}occoz}, Topological methods in modern
   mathematics (Stony Brook, NY, 1991), Publish or Perish, Houston, TX, 1993,
   pp.~467--511.

\bibitem[K]{kuratowski2}
Casimir Kuratowski, \emph{Topologie. {V}ol. {II}}, Troisi\`eme \'edition,
   corrig\`ee et compl\'et\'ee de deux appendices. Monografie Matematyczne, Tom
   21, Pa\'nstwowe Wydawnictwo Naukowe, Warsaw, 1961.


\bibitem[MSS]{mss}
Ricardo Ma{\~n}{\'e}, Paulo Sad, and Dennis Sullivan, \emph{On the dynamics of
   rational maps}, Ann. Sci. \'Ecole Norm. Sup. (4) \textbf{16} (1983), no.~2,
   193--217.

\bibitem[M]{jackdynamicsthird}
John Milnor, \emph{Dynamics in one complex variable}, third ed., Annals of
   Mathematics Studies, vol. 160, Princeton University Press, Princeton, NJ,
   2006.

\bibitem[Na]{continuumtheory}
Sam~B. Nadler, Jr., \emph{Continuum theory. {A}n introduction}, Monographs and
   Textbooks in Pure and Applied Mathematics, vol. 158, Marcel Dekker Inc., New
   York, 1992.

\bibitem[Ne]{newmanplanetopology}
M.~H.~A. Newman, \emph{Elements of the topology of plane sets of points},
   Cambridge, At the University Press, 1951, 2nd ed.

\bibitem[P]{pommerenke}
Christian Pommerenke, \emph{Boundary behaviour of conformal maps}, Grundlehren
   der Mathematischen Wissenschaften, vol. 299, Springer-Verlag, Berlin, 1992.


\bibitem[R1]{thesis}
Lasse Rempe, \emph{Dynamics of exponential maps}, doctoral thesis,
   Christian-Albrechts-Universit\"at Kiel, 2003,
   \href{http://e-diss.uni-kiel.de/diss_781}{{\tt
   http://e-diss.uni-kiel.de/diss\_781/}}.


\bibitem[R2]{topescapingnew}
\bysame, \emph{Topological dynamics of exponential maps on their escaping
   sets}, Ergodic Theory Dynam.~Systems \textbf{26} (2006), no.~6, 1939--1975,
   \mbox{\href{http://www.arXiv.org/abs/math.DS/0309107}{arXiv:math.DS/0309107}%
}.

\bibitem[R3]{localconnectivity}
\bysame, \emph{Prime ends and local connectivity}, Preprint, 2007,
   \mbox{\href{http://www.arXiv.org/abs/math.GN/0309022}{arXiv:math.GN/0309022}%
}, to appear in Bull.\ London Math.\ Soc.


\bibitem[R4]{rationalfibers}
Lasse Rempe, \emph{Triviality of rational fibers in exponential parameter
   space}, in preparation.


\bibitem[RS1]{bifurcations_new}
Lasse Rempe and Dierk Schleicher, \emph{Bifurcations in the space of
   exponential maps}, 2007,
   \mbox{\href{http://www.arXiv.org/abs/math.DS/0311480v6}{arXiv:math.DS/031148%
0v6}}, submitted for publication; previous version appeared in Stony Brook IMS
   Preprint \#2004/03.

\bibitem[RS2]{expcombinatorics}
\bysame, \emph{Combinatorics of bifurcations in exponential parameter space},
   Transcendental Dynamics and Complex Analysis (P.~Rippon and G.~Stallard,
   eds.), London Math. Soc. Lecture Note Ser. 348, Cambridge Univ. Press, 2008,
   \mbox{\href{http://www.arXiv.org/abs/math.DS/0408011}{arXiv:math.DS/0408011}%
}, pp.~317--370.

\bibitem[S1]{habil}
Dierk Schleicher, \emph{On the dynamics of iterated exponential maps},
   Habilitation thesis, TU M\"unchen, 1999.

\bibitem[S2]{expattracting}
\bysame,
   \emph{\href{http://www.math.helsinki.fi/Annales/Vol28/schleich.html}{Attract%
ing dynamics of exponential maps}}, Ann. Acad. Sci. Fenn. Math. \textbf{28}
   (2003), 3--34.

\bibitem[S3]{cras}
\bysame, \emph{Hyperbolic components in exponential parameter space}, C. R.
   Math. Acad. Sci. Paris \textbf{339} (2004), no.~3, 223--228.

\bibitem[S4]{fibers}
\bysame, \emph{On fibers and local connectivity of {M}andelbrot and {M}ultibrot
   sets}, Fractal geometry and applications: a jubilee of Beno\^\i t Mandelbrot.
   Part 1, Proc. Sympos. Pure Math., vol.~72, Amer. Math. Soc., Providence, RI,
   2004, pp.~477--517.

\bibitem[S5]{intaddrnew}
\bysame, \emph{Internal addresses in the {M}andelbrot set and irreducibility of
   polynomials}, Preprint, 2007,
   \mbox{\href{http://www.arXiv.org/abs/math.DS/9411238v2}{arXiv:math.DS/941123%
8v2}}, updated version of Stony Brook IMS Preprint \#1994/19.

\bibitem[SZ]{expescaping}
Dierk Schleicher and Johannes Zimmer, \emph{Escaping points of exponential
   maps}, J. London Math. Soc. (2) \textbf{67} (2003), no.~2, 380--400.

\bibitem[TL]{tanleiparabolics}
Tan Lei, \emph{Local properties of the {M}andelbrot set at parabolic points},
   The Mandelbrot set, theme and variations, London Math. Soc. Lecture Note
   Ser., vol. 274, Cambridge Univ. Press, Cambridge, 2000, pp.~133--160.


\bibitem[T]{thurstonlaminations}
William Thurston, \emph{The combinatorics of iterated rational maps}, Preprint,
   1985.

\end{thebibliography}

\providecommand{\href}[2]{#2}\def\polhk#1{\setbox0=\hbox{#1}{\ooalign{\hidewidth 
\lower1.5ex\hbox{`}\hidewidth\crcr\unhbox0}}}
   \def\polhk#1{\setbox0=\hbox{#1}{ \ooalign{\hidewidth
   \lower1.5ex\hbox{`}\hidewidth\crcr\unhbox0}}} \input{cyracc.def} \def\j{{\u
   i}} \def\J{{\u I}} \newfont{\cyrit}{wncyi10 at 12pt}\def\cprime{$'$}
\providecommand{\bysame}{\leavevmode\hbox to3em{\hrulefill}\thinspace}

\end{document}